\documentclass{article}
\usepackage{graphicx,natbib}
\usepackage{setspace}

\usepackage{verbatim}
\usepackage{geometry}          
\usepackage[parfill]{parskip}     
\usepackage{subfigure}
\usepackage{amssymb,amsmath}
\usepackage{epstopdf}
\usepackage{graphicx}
\usepackage{setspace}
\usepackage{sidecap}
\newcommand{\ds}[1]{\displaystyle#1}

\begin{document}
\title{A Mechanism For Stabilization of Dynamics in Nonlinear Systems with Different time Scales}

\author{Raquel M. L\'opez$^{1,}$
\footnote{ \texttt{{\small rlopez14@asu.edu, $^\dag$sks@asu.edu, $^\ddag$erika.camacho@asu.edu} } } \\
  Sergei K. Suslov$^{2,\,\,\dag}$\\
  Erika T. Camacho$^{3,\,\,\ddag}$ \\ 
 {\small $^1$ Mathematical, Computational \& Modeling Sciences Center,}\\
 {\small School of Mathematical and Statistical Sciences,} 
 {\small Arizona State University, Tempe, AZ}\\
 {\small $^2$ School of Mathematical and Statistical Sciences,}\\
 {\small Arizona State University, Tempe, AZ}\\ 
 {\small $^3$ Division of Mathematics \& Natural Sciences,} \\ 
 {\small New College of Interdisciplinary Arts \& Sciences, Arizona State University, Phoenix, AZ} }

\date{\today}

\maketitle

\begin{abstract}
There are many natural, physical, and biological systems that exhibit multiple time scales.  For example, the dynamics of a population of ticks can be described in continuous time during their individual life cycle yet discrete time is used to describe the generation of offspring.  These characteristics cause the population levels to be reset periodically.  A similar phenomenon can be observed in a sociological college drinking model in which the population is reset by the incoming class each year, as described in the 2006 work of Camacho et al.  With the latter as our motivation we analytically and numerically investigate the mechanism by which solutions in certain systems with this resetting conditions stabilize.  We further utilize the sociological college drinking model as an analogue to analyze certain one-dimensional and two-dimensional nonlinear systems, as we attempt to generalize our results to higher dimensions. \\

\end{abstract}

\section*{Introduction}
In many dynamical systems, especially systems found in the biological and natural sciences, it is common to find  multiple time scales embedded in the same system all of which might be continuous or include a combination of both discrete and continuous time.  Dynamical systems with multiple time scales are also called singularly perturbed systems and slow-fast systems.  An excellent example of a slow-fast type is the forced Van der Pol equation [5], and is the first system in which the phenomenon of chaos was first observed. However for our purposes and interests, we restrict only to discrete--continuous time scaled systems.  Examples of these arise mostly in ecology.  To illustrate, recall the famous system of equations known as the \emph{Lotka--Volterra model} [ 1, 3 ].  In this model, $x(t)$ and $y(t)$ denote, the prey and predator populations, respectively, at any time $t$.  Now let us imagine the following prey-predator relationship: The predators are sharks, and the prey are small fish. If there were no more small fish, one might expect that the population of sharks would decline in quantity, according to the differential equation 
\begin{equation}
\begin{array}{rcl}
\ds\frac{dy}{dt} & = & -cy, \qquad c > 0.  
\label{in1}
\end{array}
\end{equation}
When small fish are present in the habitat or environment, however, it seems more reasonable that the number of encounters between these two species per unit time is jointly proportional to their populations $x$ and $y$, or in other words, proportional to the product $xy$, by the law of mass action [4].  The addition of this last rate yields an equation for the shark population, and so in a similar fashion, the dynamics of the prey population $\ds\frac{dx}{dt}$ is derived, giving us the nonlinear system of differential equations known as the prey--predator model (2), where $a$, $b$, $c$ and $d$ are all positive constants.  
\begin{equation}
\begin{array}{rcl}
\ds\frac{dx}{dt} & = & ax - bxy\\[.1in] 
\ds\frac{dy}{dt} & = & -cy + dxy 
\label{in2}
\end{array}
\end{equation}
Although the explicit solution of system (2) is not known, it can be analyzed quantitatively and qualitatively [3], that is we can establish some information about the behavior of its solutions, by seeking a relation for $x$ and $y$ ignoring temporal space ( a better explanation of this can be found in [1] ).  Moreover, we can look at this model in a slightly different way by making the growth rate of the prey time dependent, in other words, not constant.  If we change the prey population slightly as follows, 
\begin{equation}
\begin{array}{rcl}
\ds\frac{dx}{dt} & = & a(t)x -bxy,
\label{in3}
\end{array}
\end{equation}
then the population dynamics still evolve continuously in time, but what differs now is the birthrate of the prey, denoted by some given $a(t)$ which evolves discretely in time.  This is due to normal ecological changes in temperature or the environment, and as a result the population behaves periodically according to this birthrate.  In addition, we also find these type of phenomena in fishing models with harvesting such as the one found in [1] given by the generalized equation
\begin{equation}
\begin{array}{rcl}
\ds\frac{dx}{dt}& = & f(x)- H(t).
\label{in4}
\end{array} 
\end{equation}
Here $H(t)$ is a non--constant rate of members ( fish ) harvested per unit time, and as in the \emph{Lotka--Volterra} model, the population of fish evolves continuously, yet due to the harvesting rate $H(t)$, the population resets periodically, and this periodic reseting is dependent, more realistically, on harvesting seasons.  \\
As observed, such examples of multi-scaled systems are naturally found in nature, but such systems are also found in sociological models like a deterministic Model of College Drinking Patterns analyzed by Camacho, et al. (2006) [6] ( See also [10] and references therein ).  This model has become the onset point of this prospectus journey.  A brief description of the continuous deterministic model given by equations (5)--(8) below, is as follows.  First a four equation system is given, inclusive of an embedded discrete difference equation $R_j(t_{k+1})$ (9) which is originally called the resetting equation due to Scribner, et al. [10].  This equation is crucial to describe periodicity and resetting of a population as it takes care of replacing the total number of students that left at the end of each academic year, back into the system at the beginning of each academic year by keeping the total population in each class constant.  Through this equation, every year current variable values are reset, resulting in a non--smooth set of solution curves.  In other words, the population in each class ( $N_1, N_2, N_3, N_4$ ) gets reset to bring the total population down to the initial data $N(0)$.  This resetting force leads to stabilization but at a different level than the equilibria from the continuous system absent from resetting forces. \\
\begin{eqnarray}
\ds\frac{dN_1}{dt} &= & -d_1N_1 + r_{21}N_2 + r_{31}N_3 - s_{12} \ds\frac{N_1 N_2}{N}- n_{12}\ds\frac{N_1N_2}{N}\\[.15in] 
\nonumber\ds\frac{dN_2}{dt}  &=&  -d_2 N_2 - r_{21}N_2 - r_{23} N_2 - r_{24}N_2 + r_{42}N_4 + s_{12}\ds\frac{N_1N_2}{N}\\[.1in] 
& & -s_{23}\ds\frac{N_2N_3}{N} + (s_{42} - s_{24} )\ds\frac{N_2N_4}{N} + n_{12}\ds\frac{N_1N_2}{N} - n_{24} \ds\frac{N_2 N_4}{N}\\[.15in]
\ds\frac{dN_3}{dt}& =&  -d_3N_3 + r_{23}N_2 - r_{31}N_3 + r_{43}N_4 + s_{23}\ds\frac{N_2N_3}{N} + s_{43}\ds\frac{N_4N_3}{N}\\[.15in]
\nonumber\ds\frac{dN_4}{dt}  &=&  -d_4N_4 + r_{24}N_2 - r_{42}N_4 - r_{43}N_4 + (s_{24} - s_{42} )\ds\frac{N_2N_4}{N}\\[.1in] 
&&-s_{43}\ds\frac{N_4N_3}{N} + n_{24}\ds\frac{N_2N_4}{N}
\label{in5}
\end{eqnarray}
\begin{eqnarray}
R_j (t_{k+1} )  & = &  c_j [N(0)-  \sum_{i=1}^4 N_i(t^-_k)]
\label{in6}
\end{eqnarray}
\begin{eqnarray}
N(t) &  = & N_1 (t) + N_2 (t) + N_3 (t) + N_4 (t)
\label{in7}
\end{eqnarray}
From the resetting equation, $t_{k+1}$ and $t^-_k$ represent the beginning of year $k+1$ and the end of year $k$, respectively, and the negative superscript indicates that $t_{k+1}$ is dependent on the previous initial time $t_k$.  It can also be noted that $N(0)$ is the total initial population, and  $c_j  =  \ds\frac{N{_j}(0)}{N(0)}$ for $j = \{1, 2, 3, 4\}$ is a percentage of the total population which enters each class.  Moreover, in the modified model discussed by Camacho et al. [6], numerical tools were the only option available to analyze the nonlinearity of the system, hence numerical simulations have shown that stabilization of the solutions occurred even when the continuous system without resetting did not stabilize ( i.e. solutions of the system with resetting present did not attain the stable equilibrium belonging to the continuous system absent of resetting ).  It was also noticed that for certain initial conditions the system ( equations (5)--(8) ) stabilized at some artificial equilibrium.  Analytically, it was almost impossible to investigate the behavior in the system as a whole, and the worse was finding an explicit solution.  This facts have motivated us to search for analytic mechanisms which lead solutions of nonlinear systems with different time scales to stabilize periodically in time.  In this quest, many questions have arised.  $(i)$ Under what conditions will the system stabilize and $(ii)$ is there a way of proving analytically that stabilization of the system occurs?  In the latter case, due to the nonlinearities of the four dimensional alcohol model, proving and obtaining conditions for stabilization analytically are not easy tasks.  Thus in this project we have focused on searching for mechanisms which stabilize system dynamics, from the simplest model available moving on to more complicated systems of larger dimensions once some conclusions about simpler models are established.  Thus, we begin investigating this mechanisms in nonlinear models case by case. 
  
\section*{Analytical Investigation of Simplest Nonlinear Models \\(Summary of Results)}
In order to describe a mechanism of stabilization in nonlinear dynamical systems, we need to start from the simplest linear and  nonlinear models available.  In this prospectus summary we have considered elementary nonlinear cases both one--dimensional and two--dimensional.  With the following introductory but significant examples, we show that specific resetting conditions and parameter values give rise to the stabilization of these systems. The simplest one--dimensional available case considered are exponential growth and decay (Malthus Law) where  
\begin{eqnarray*}
\ds\frac{dx}{dt} &=& \alpha x \quad \mbox{ and } \quad  \ds\frac{dx}{dt} = -\alpha x
\label{in8}
\end{eqnarray*}
respectively, are the well known differential equations describing exponential growth and decay, with their very well known solutions namely 
\begin{eqnarray*}
x(t) &=& x(0)e^{\alpha t} \quad \mbox{ and } \quad  x(t) = x(0)e^{-\alpha t}.
\label{in9}
\end{eqnarray*}
Secondly we have considered both, a logistic model
\begin{eqnarray}
 \ds\frac{dx}{dt} & = & \alpha x \left(1-\ds\frac{x}{\beta} \right)
 \label{in10}
 \end{eqnarray}
and a Gompertz (1825) model (taken from [1])
\begin{eqnarray}
\ds\frac{dx}{dt} &=& -\alpha x \mbox{ln}\left(x\right).
 \label{in11}
\end{eqnarray}
For these models we have artificially created a resetting condition given by
\begin{eqnarray}
x\left(T_{n+1}\right) &=& \gamma x\left(T^-_{n+1}\right)
\label{in12}
  \end{eqnarray}
where $\gamma$ is a resetting parameter, and $0< \gamma <1$.
On the other hand, in the two--dimensional case we have considered a Logistic Coupled system [3]
\begin{eqnarray}
\ds\frac{dx}{dt} & = &\alpha x \left(1-\ds\frac{x}{\beta}\right) \\[.1in] 
\ds\frac{dy}{dt} & = & \beta xy
\label{in13}
\end{eqnarray}
and a Gompertz-type System [1] 
\begin{eqnarray}
\ds\frac{dx}{dt} & = &-\alpha x \mbox{ln}\left(x\right) \\[.2in] 
\ds\frac{dy}{dt} & = & -\beta y \mbox{ln} \left(x\right)
\label{in14}
\end{eqnarray}
for which we have chosen the resetting condition
\begin{eqnarray}
\left[\begin{matrix} 
 x(T_{n +1}) \\
 y(T_{n+1})
\end{matrix} \right]
&=&
\left(\begin{matrix} 
\lambda & \xi \\ 
\gamma & \mu 
\end{matrix} \right)\cdot
 \left[\begin{matrix} 
 x(T^-_{n+1}) \\ 
 y(T^-_{n+1})
\end{matrix} \right]
\label{in15}
\end{eqnarray}
(A more detailed explanation of the mechanisms chosen is given in the proceeding parts of this prospectus summary).
\subsection*{Summary}
Our illustrations have shown that the  resetting mechanisms we investigate, certainly stabilize nonlinear system dynamics.  As a hypothesis we propose the generalization, that any nonlinear system can be analytically investigated by following simply the theory of dynamical systems along with an advanced calculus approach as follows.  It has been established that for a nonlinear autonomous system of ordinary differential equations 
\begin{eqnarray*}
{\bf \dot x} &= &{\bf f(x)}
\end{eqnarray*}
which has a given initial condition,
\begin{eqnarray*}
{\bf x(0)} & = & {\bf x_0} 
\end{eqnarray*}
the fundamental \emph{ Existence--Uniqueness Theorem }[4] states that as long as there exists an open subset $S$ of ${\bf R}^n$, which contains ${\bf x_0}$ along with the requirement ${\bf f} \in C^1(S)$, then there must be some $a>0$ such that the initial value problem (IVP) stated
has locally an explicit analytical, unique solution ${\bf x(t)}$ on some interval $[ -a, a ]$.  In addition to the uniqueness of the solution of the system, resetting implies intersection of periodic orbits in the state space of the continuous dynamical system, with certain lower dimensional subspace called Poincar\'e Sections, transverse to the flow of the system.  Initial data are thus periodic orbits with initial conditions on the Poincar\'e section thereby formulating a Poincar\'e map 
\begin{eqnarray*}
{\bf x_{n+1}} &=& {\bf f(x_n)}. 
\end{eqnarray*}
This iteration process or map allows us to further, invoke the contraction mapping theorem, and thus we can claim that the dynamics of the system are guaranteed to stabilize, when a generated sequence ( Poincar\'e Sections ) of initial data converges to a fixed point of the corresponding discrete mapping.  But in general, for $n$-dimensional systems, it is very difficult to find a domain of initial data, where the discrete map might converge analytically.  Investigation of the \lq\lq size\rq\rq of this region is a typical inverse problem which will be investigated numerically in the following sections.  \\
Furthermore, there are many possibilities and applications to this work.  Higher dimensional systems for example, could be analyzed with monotone systems theory, by means of explaining the behavior of high dimensional systems by that of one dimensional flow.  But in general we want to find how the interplay of initial conditions and resetting definitions affect the behavior of solutions in the long run.  Ultimately we want to find the set of initial conditions that leads to stabilization of dynamics by means of a convergent discrete mapping or Poincar\'e Mapping.  Our future considerations include the investigation of any type of model in which resetting occurs and ideally we want to find applications to ecological and sociological models.  It is also of potential interest to consider models in which the periodic administration of a drug, such as chemotherapy for instance, is given to cancer patients.  \\\\


\noindent \textbf{Acknowledgments.\/} This paper is written as a part of the
summer 2008 program on analysis of the Mathematical and Theoretical Biology
Institute (MTBI) at Arizona State University, now in conjunction with the Mathematical, Computational 
\& Modeling Sciences Center (MCMSC). The MTBI/SUMS Summer
Undergraduate Research Program is supported by The National Science
Foundation (DMS-0502349), The National Security Agency (DOD-H982300710096),
The Sloan Foundation, and Arizona State University. The authors are grateful
to Professor Carlos Castillo-Ch\'avez for support, reference and valuable comments.
\clearpage 
{\center {\bf \Large Appendix: Examples of Numerical Results}
\begin{figure}[htbp]
\centering
\includegraphics[width=6in, height = 3.5in]{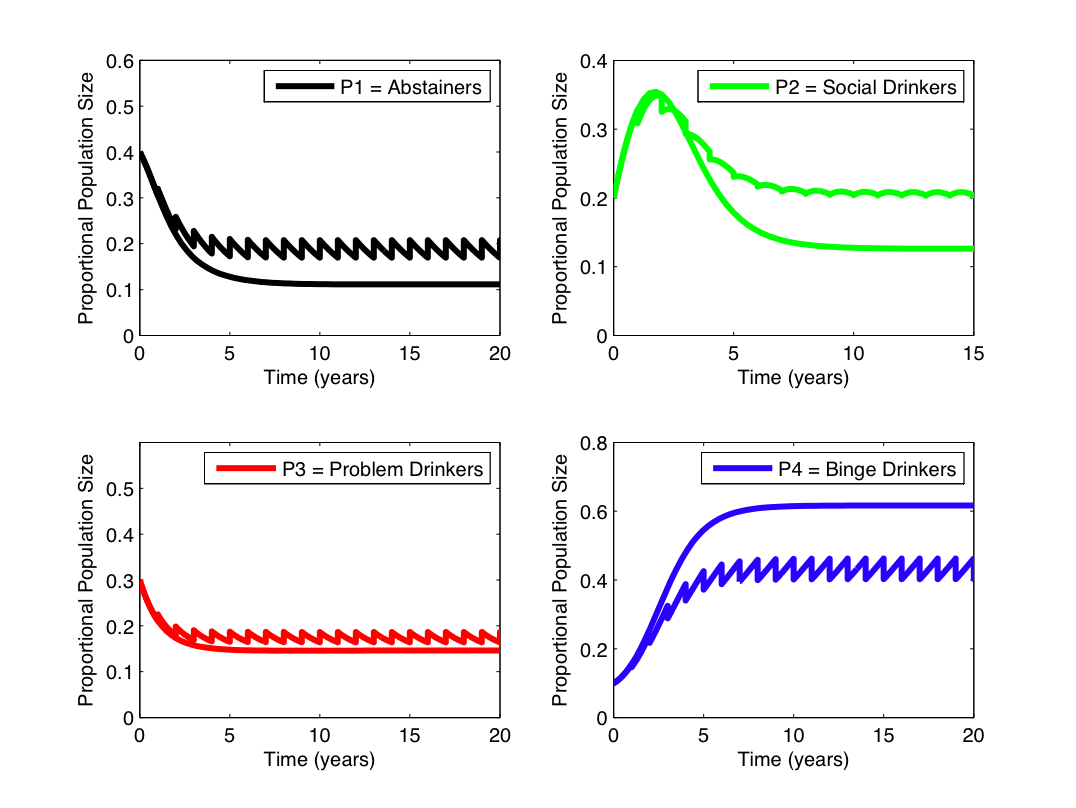} 
  \caption{{\bf College Drinking Model}.  In this figure, the smooth curves represent solutions of the continuous system ( equations (5)--(8) ) without resetting and the the non-smooth curve represent the solution of the system with resetting.  Note that the stabilization of the continuous system with resetting occurs at a level different from the equilibria of the continuous system without resetting.  That is, when resetting is introduced the system tends to stabilize to an equilibria of its own.}
  \label{fig1}
 \end{figure}
\begin{figure}[htbp]
\centering
\begin{tabular}{ccc}
\includegraphics[width=0.45\linewidth]{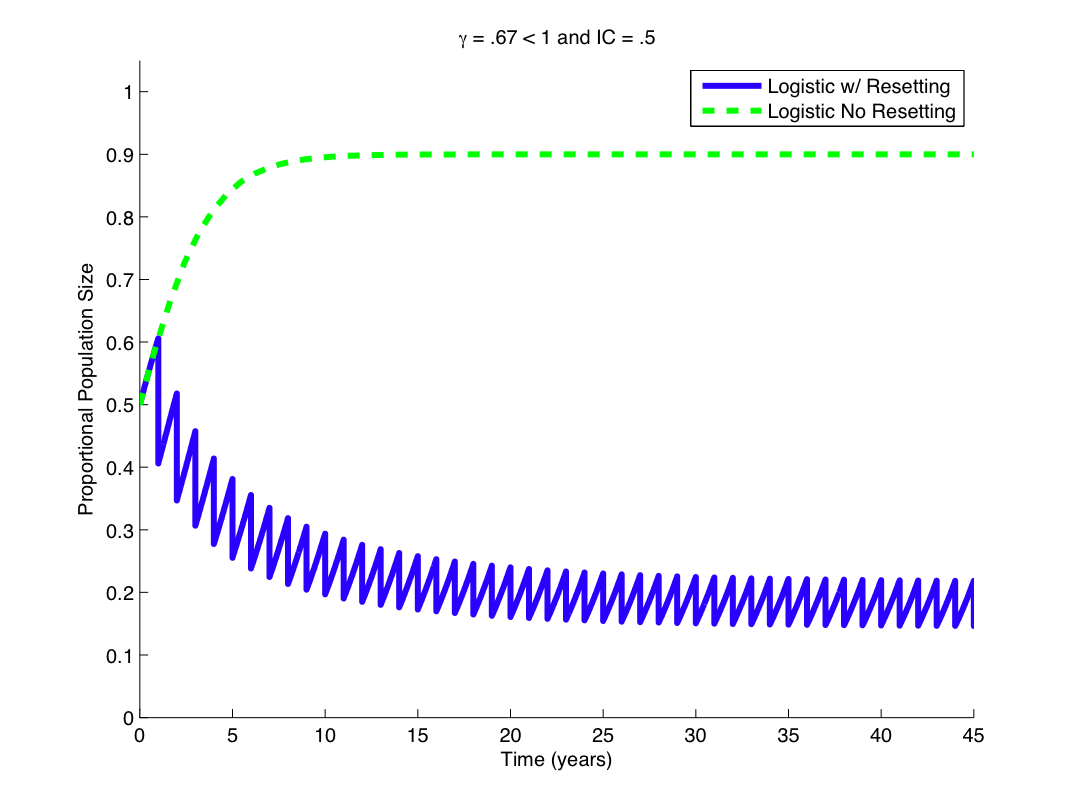} & & 
\includegraphics[width=0.45\linewidth]{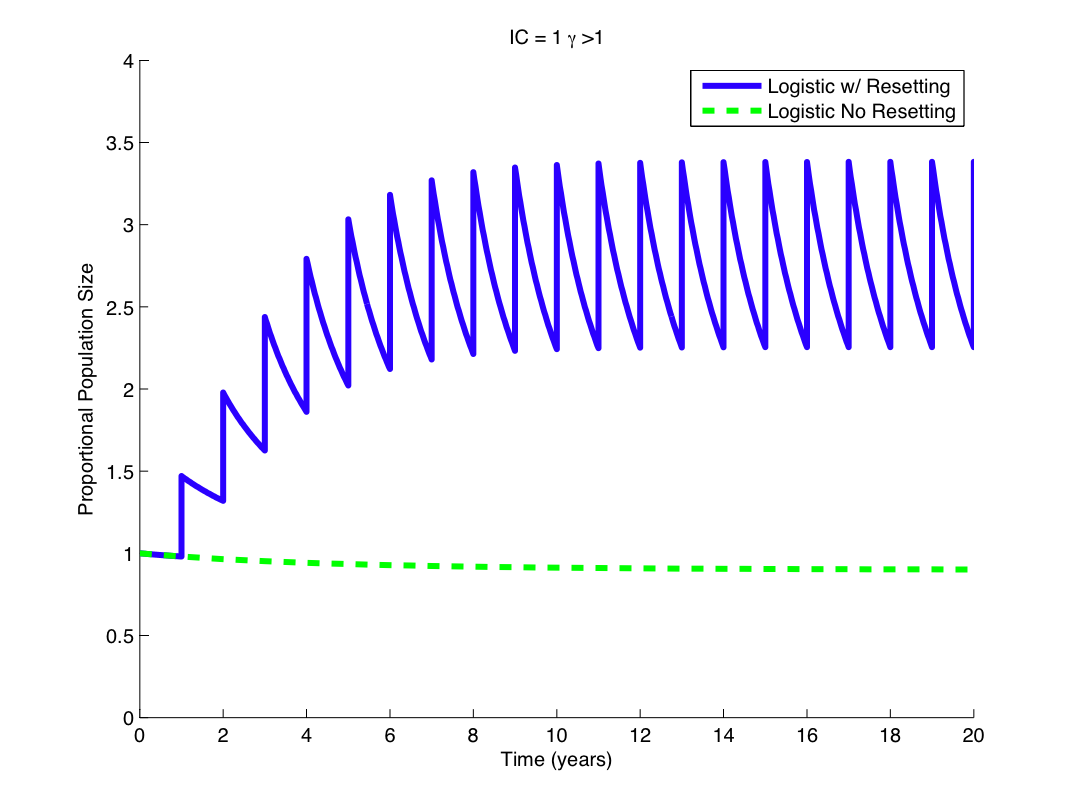}
  \end{tabular}
     \caption{{\bf One--Dimensional Logistic Model with Resetting}  Using the resetting condition chosen for the one--dimensional logistic model ( equation (11) ), numerical simulations show the stabilization expected, which match those observed by Camacho et al.  Here again, the smooth curves correspond to the continuous system without resetting.  At some time, the resetting rule comes into effect, and the solutions become nonsmooth.  Notice in the first case for a resetting parameter $\gamma = 0.67 < 1$ and initial condition $x_0 = 0.5$ ( $0<x_0<0.9$, where $x^*_{stable} = 0.9$, and $x^*_{unst} = 0$), the dynamics of the system with resetting stabilize between the equilibrium of the continuous system, yet they follow an equilibrium level of its own.  In the second case for initial conditions above stable equilibrium, we see the same phenomena ocurring.  It is very interesting to note that through this resetting mechanism ( or rule ) we have chosen for the simple one--dimensional model, stabilization is not the level of the system without resetting. }
\label{fig2}
\end{figure}
\begin{figure}[htbp]
\begin{center}
     \includegraphics[width=5.5in, height = 2in]{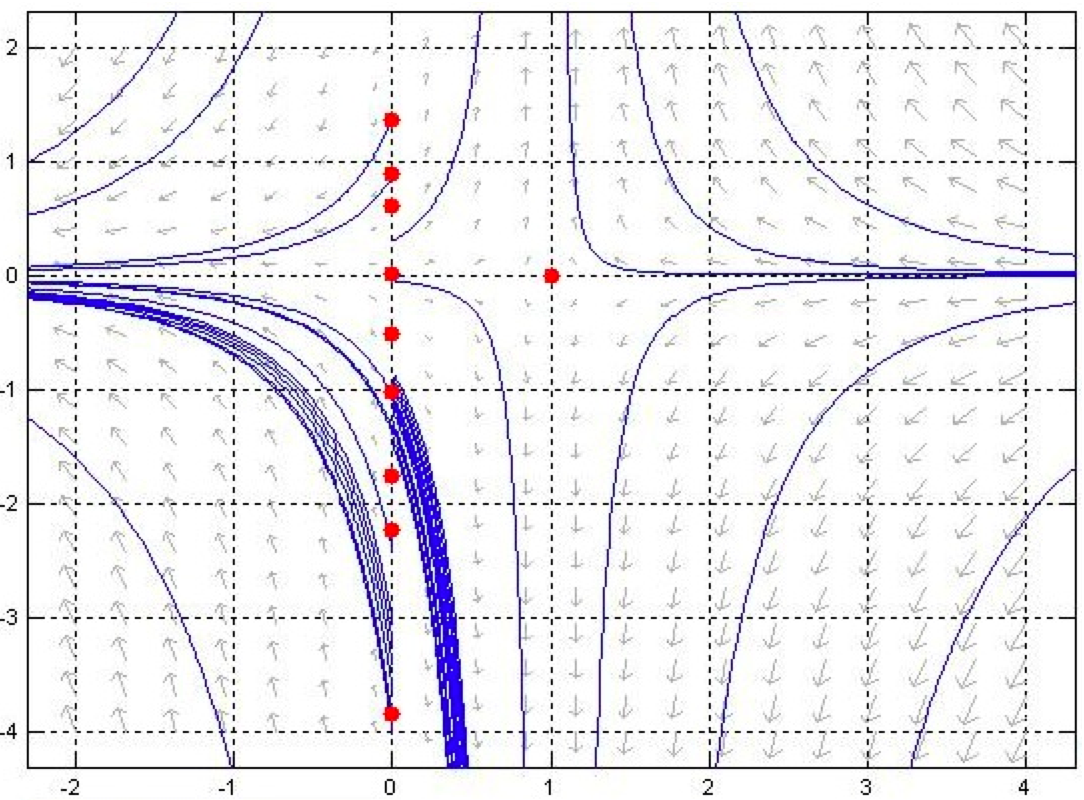}   
     \label{fig3}
     \caption{\small {\bf Phase Portrait of Logistic Coupled system absent from resetting}  The chosen coupled system $\dot x = \alpha x (1-x/\beta)$,  $\dot y = \beta xy$ for $0<\alpha$,
      $\beta<1$ which has explicit continuous solutions given by $x(t) = \ds\frac{\beta x(0)}{x(0) + (\beta - x(0))e^{-\alpha t}}$ \quad and \quad $y(t) = y(0)(1+x(0)(e^{\alpha t} -1)^{\beta/\alpha} )$ has a phase portrait where equilibria $(x^*, y^*) = (1, 0)$, and  $(0, y^*)$ are a saddle point and a set of degenerate isolated unstable equilibria, respectively.  Every trajectory is unstable, as the behavior of the system without resetting exhibits a saddle node, unless the initial conditions fall on the stable manifold $(1, 0)$.} 
 \end{center}
 \label{fig4}
\end{figure}
     \begin{figure}[htbp]
\centering
\begin{tabular}{ccc}
     \includegraphics[width=0.45\linewidth]{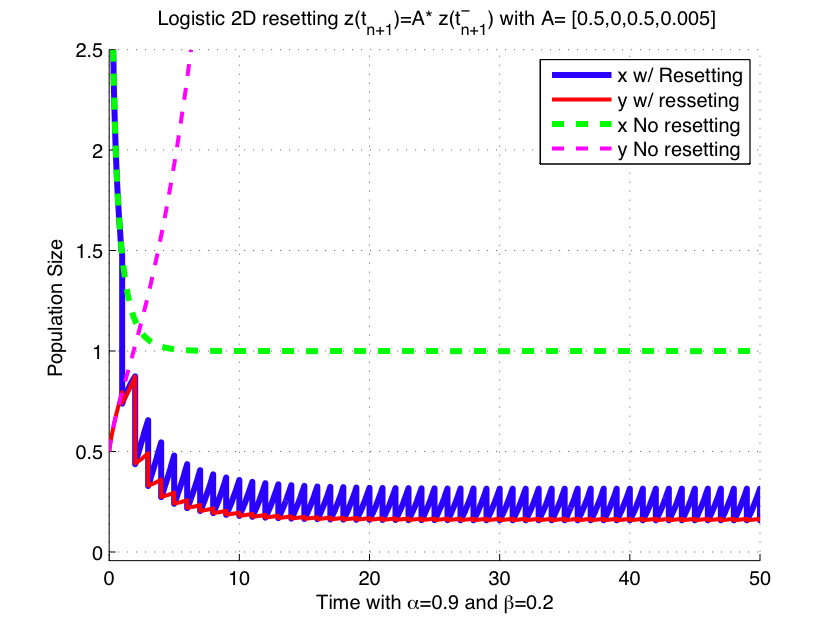} & &
     \includegraphics[width=0.45\linewidth]{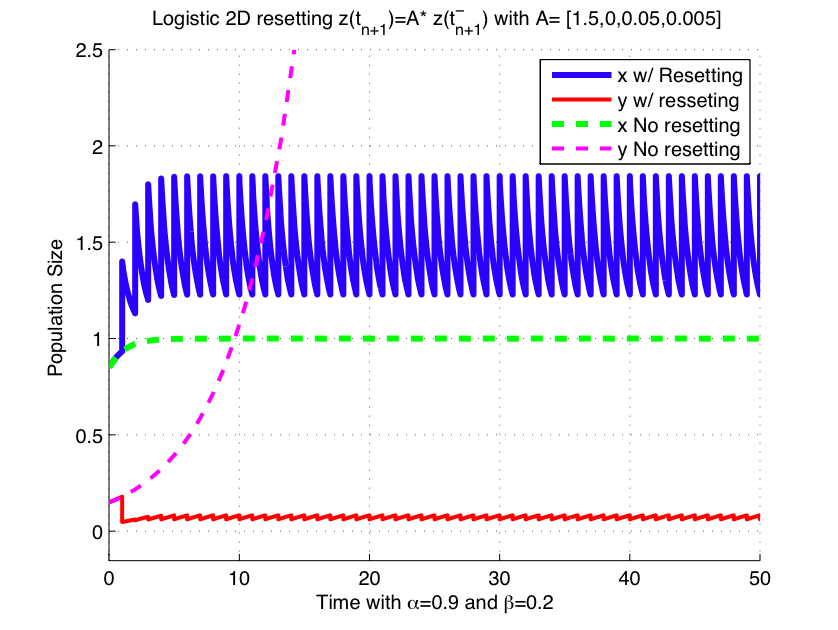}
  \end{tabular}
  \label{fig5}
     \caption{{\bf Two--Dimensional Logistic Coupled System} Simulations belonging to the two--dimensional system ( equations 14--15 ) again show, with the chosen resetting rule once again the system stabilizes, as expected for certain initial conditions and resetting parameter $\gamma$.  The solution $y(t)$ in the continuous system without resetting ( represented by the pink dashed curve ), is tending to infinity, while with resetting ( represented by the red solid curve) $y(t)$ tends to some nonzero equilibrium ( very close to zero ) but in an oscillatory manner with very small amplitude.  The solution $x(t)$ ( represented by the blue solid curve ) in the system with resetting also stabilizes to a different level than the continuous solution $x(t)$ in the system without resetting.  Notice the amplitude of the oscillations of the solutions of $x(t)$ surpass those of $y(t)$ ( the red curve ) and this is due to the choice in size of the parameters.  It is quite surprising that although the continuous system has no stable solution the introduction of a periodic resetting still causes the system to stabilize.}
     \end{figure}
     \clearpage 

 \end{document}